\documentclass[12pt]{article}
\usepackage{a4}
\usepackage{amsthm}
\usepackage{amsmath}
\usepackage{amssymb}
\usepackage{amsfonts}
\usepackage{cite}
\usepackage{epsfig}
\newtheorem{theorem}{Theorem}

\begin{document}
\title{Hadwiger meets Cayley}
\author{Jacob W. Cooper\thanks{Faculty of Informatics, Masaryk University, Botanick\'a 68A, 602 00 Brno, Czech Republic. E-mails: {\tt jcooper@mail.muni.cz}, {\tt kabela@fi.muni.cz}, {\tt dkral@fi.muni.cz} and {\tt pierron@fi.muni.cz}. All authors were supported by the MUNI Award in Science and Humanities of the Grant Agency of Masaryk university.}\and
\newcounter{lth}
\setcounter{lth}{1}
        Adam Kabela$^\fnsymbol{lth}$\and
        Daniel Kr\'al'$^\fnsymbol{lth}$\thanks{Mathematics Institute, DIMAP and Department of Computer Science, University of Warwick, Coventry CV4 7AL, UK.}\and
	Th\'eo Pierron$^\fnsymbol{lth}$}
\date{}
\maketitle

\begin{abstract}
We show that every connected $k$-chromatic graph contains at least $k^{k-2}$ spanning trees.
\end{abstract}

The present note is motivated by a corollary of the well-known conjecture of Hadwiger on $k$-chromatic graphs.
Recall that a graph is $k$-chromatic if it admits a proper vertex coloring using $k$ colors but not $k-1$ colors.
In 1943, Hadwiger~\cite{Had43} conjectured that every $k$-chromatic graph contains $K_k$ as a minor.
The conjecture
remains wide open despite the effort of many researchers,
and is considered to be one of the deepest unsolved problems in graph theory,
see e.g.~\cite{BohKT11,BolCE80,Dir52,Fox10,Kaw07,KawT05,Kos84,RobST93,Wag37}
and the surveys~\cite{Sey16,Tof96}.
The conjecture has been confirmed for $k\le 6$~\cite{RobST93}.
If true, the conjecture would imply that
every connected $k$-chromatic graph contains at least $k^{k-2}$ spanning trees
(by Cayley's formula~\cite{Bor61,Cay89} on the number of spanning trees of complete graphs).
We give a direct and short proof of this statement.
This answers a problem recently posted by Sivaraman~\cite[Problem 8]{Siv20}.

\begin{theorem}
\label{thm:main}
Every connected $k$-chromatic graph contains at least $k^{k-2}$ spanning trees.
\end{theorem}

Before proving the result, we fix some notation.
Given a graph $G$, we denote by $G-v$ the graph obtained from $G$ by removing vertex $v$, and
by $G/v_1\ldots v_\ell$ the graph obtained from $G$ by identifying vertices $v_1,\ldots,v_{\ell}$ and removing all newly created loops and parallel edges so as to yield a simple graph.

\begin{proof}
We proceed by induction on the number of vertices.
Clearly, the statement is satisfied if $k=1$ or $k=2$, which covers the base of the induction.
We assume that $k\ge 3$ and fix a connected $k$-chromatic graph $G$.
We first show that $G$ can be assumed to be $2$-connected, have minimum degree at least $k-1$, and contain no $K_k$ as a subgraph.

{\bf The graph $G$ is 2-connected.}
Indeed, suppose that $G$ contains a vertex $v$ for which $G-v$ is disconnected.
For every connected component of $G - v$,
consider the subgraph of $G$ induced by $v$ and the vertices of the component.
Note that at least one of the considered subgraphs has the same chromatic number as $G$,
and hence this subgraph contains at least $k^{k-2}$ spanning trees by the induction hypothesis.
The statement follows by observing that the number of spanning trees of $G$
is equal to the product of the numbers of spanning trees of the considered subgraphs.

{\bf The minimum degree of $G$ is at least $k-1$.}
Suppose that $G$ contains a vertex $v$ of degree at most $k-2$.
Note that the graph $G-v$ is connected (since $G$ is $2$-connected) and
its chromatic number is at least $k$ (since $G$ cannot be colored using $k-1$ colors).
By the induction hypothesis, $G-v$ contains at least $k^{k-2}$ spanning trees.
Each of these trees can be extended to a spanning tree of $G$
by adding an arbitrary edge incident with $v$, and
the obtained spanning trees are all distinct.

{\bf The graph $G$ does not contain $K_k$.}
Suppose that $G$ contains $K_k$ as a subgraph (possibly $G$ is $K_k$).
By Cayley's formula, this subgraph contains $k^{k-2}$ spanning trees.
The statement follows by noting that each of these trees can be trivially extended to a spanning tree of $G$ and
the resulting spanning trees are all distinct.

We next choose an arbitrary vertex of $G$, say $v$.
Let $d$ be the degree of $v$ and $F$ be the set of all edges incident with $v$;
in particular, $|F| = d$ and $d \ge k-1$.
By assumption, the graph $G - v$ is connected, and its chromatic number is at least $k-1$.
Hence, $G-v$ contains at least $(k-1)^{k-3}$ spanning trees by the induction hypothesis, and
each of these trees can be extended to a spanning tree of $G$ by adding an arbitrary edge of $F$.
In this way, we obtain at least $d(k-1)^{k-3}$ distinct spanning trees of $G$.
We proceed by finding additional spanning trees of $G$ (each containing at least two edges of $F$, and
therefore distinct from those derived above).
We distinguish two cases based on the value of $d$.

We first consider the case $d=k-1$.
Since $G$ does not contain $K_k$ as a subgraph,
the vertex $v$ has two non-adjacent neighbors, say $a$ and $b$.
Note that the chromatic number of the graph $(G-v)/ab$ is at least $k$;
otherwise, $G$ could be colored using $k-1$ colors by taking a coloring of $(G-v)/ab$ with $k-1$ colors and modifying it as follows.
Assign the vertices $a$ and $b$ the color of the vertex obtained by their identification, and
assign $v$ a color not appearing among its neighbors' colors (there are at most $k-2$ such colors).
Moreover, $(G-v)/ab$ is connected (as $G-v$ is connected),
and thus contains at least $k^{k-2}$ spanning trees by the induction hypothesis.
Every such tree of $(G-v)/ab$ corresponds to a spanning forest of $G-v$ with two components,
one containing $a$ and the other containing $b$.
Each of these forests can each be extended to a spanning tree of $G$
either by adding both edges $av$ and $bv$, or 
by adding an arbitrary edge from the set $F\setminus\{av,bv\}$ and
adding either the edge $av$ or the edge $bv$ (the choice between $av$ and $bv$ is determined so that
the resulting graph is a tree).
This gives at least $(d-1)k^{k-2}$ distinct additional spanning trees.
We conclude that $G$ contains at least $d(k-1)^{k-3}+(d-1)k^{k-2}$ spanning trees in total,
which is at least $k^{k-2}$ as $d=k-1$.

We next discuss the case $d\ge k$.
We show that the neighborhood of $v$ can be assumed to contain two disjoint pairs of non-adjacent vertices.
Suppose otherwise.
Since $G$ does not contain $K_k$ as a subgraph,
the neighborhood of $v$ does not contain a complete graph on $k-1$ vertices.
Hence, $d=k$ and there is a triple of pairwise non-adjacent neighbors of~$v$, say $a$, $b$ and $c$.
Following the same line of arguments as presented in the previous case,
we derive that the graph $(G-v)/abc$ is connected and its chromatic number is at least $k$.
By the induction hypothesis, $(G-v)/abc$ has at least $k^{k-2}$ spanning trees.
Each of these trees can be extended to a spanning tree of $G$
either by adding the three edges $av$, $bv$ and $cv$, or
by adding any edge from the set $F\setminus\{av,bv,cv\}$ and two edges among $av$, $bv$ and $cv$ so that
the resulting graph is a tree.
This yields at least $(d-2)k^{k-2}$ distinct additional spanning trees of $G$, and
thus $G$ contains at least 
$k(k-1)^{k-3}+(k-2)k^{k-2}$ spanning trees in total,
which is at least $k^{k-2}$ as desired.

Hence, we can assume that $v$ has two disjoint pairs of non-adjacent neighbors, say $a,b$, and $c,d$.
Since the chromatic number of $G-v$ is at least $k-1$,
the chromatic number of each of the connected graphs $(G-v)/ab$ and $(G-v)/cd$ is also at least $k-1$.
By the induction hypothesis, each of these two graphs contains at least $(k-1)^{k-3}$ spanning trees.
Every spanning tree of $(G-v)/ab$ can be extended to a spanning tree of $G$
either by adding both edges $av$ and $bv$, or
by adding any edge from the set $F\setminus\{av,bv\}$ and adding either $av$ or $bv$.
Similarly, each spanning tree of $(G-v)/cd$ can be extended
either by adding $cv$ and $dv$, or
by adding any edge from the set $F\setminus\{av,bv,cv,dv\}$ and adding either $cv$ or $dv$.
In this way, we obtain at least $(d-1)(k-1)^{k-3} + (d-3)(k-1)^{k-3}$ distinct spanning trees of $G$ (those obtained
from the spanning trees of $(G-v)/ab$ contain the edge $av$ or $bv$ but the latter do not).
This together with the $d(k-1)^{k-3}$ initially constructed spanning trees implies that
$G$ contains at least
\[(3d-4)(k-1)^{k-3}\ge (3k-4)(k-1)^{k-3}\ge k^{k-2}\]
spanning trees, where the last inequality holds as $k\ge 3$.
\end{proof}

We remark that a careful inspection of the proof yields that
the only $k$-chromatic graphs with exactly $k^{k-2}$ spanning trees
are those obtained from $K_k$ by iteratively adding pendant edges.

\bibliographystyle{bibstyle}
\bibliography{cayleychi}

\end{document}